\documentclass[12pt]{amsart}

\usepackage{latexsym}
\usepackage{amsfonts}
\usepackage{amsthm}
\usepackage{graphicx}
\usepackage{amssymb}
\usepackage{amsmath}
\usepackage{amssymb}
\usepackage{color}

\usepackage{hyperref}
\usepackage[top=1in, bottom=1.1in, left=1in, right=1in]{geometry}

\newtheorem{conjecture}{Conjecture}

\newcommand{\Z}{\mathbb{Z}}

\newcommand{\Q}{\mathbb{Q}}

\author[]{Vladimir Shpilrain}
\address{Department of Mathematics, The City  College  of New York, New York,
NY 10031} \email{shpilrain@yahoo.com}

\begin{document}

\title[Knot invariants]{Knot invariants from representations of braids\\
 by automorphisms of a free group}

\begin{abstract}
We describe an alternative way of computing Alexander polynomials of knots/links, based on the Artin representation of the corresponding braids by automorphisms of a free group.
Then we apply the same method to other representations of braid groups discovered by Wada and compare the corresponding isotopic invariants to Alexander polynomials.

\end{abstract}

\maketitle

\hfill{\small \it In memory of Vitaly Romankov}

\section{Introduction}

Our procedure for obtaining isotopic invariants of a link is based on the following well-known facts; a general reference is the monograph \cite{Birman}.
\medskip

$\bullet$  Every link is a closed braid. Denote by $B_n$ the braid group on $n$ strands.
\medskip

$\bullet$ Two braids $\beta_1 \in B_n$ and $\beta_2 \in B_m$ produce isotopically equivalent links when closed if and only if the braid word $\beta_1$ can be taken to the braid word $\beta_2$ by a sequence of Markov moves. (See our Section \ref{Markov} for more details.)
\medskip

$\bullet$ The abelianization (i.e., the factor group by the commutator subgroup) of any braid group $B_n$ is an infinite cyclic group.
\medskip

$\bullet$  Every braid group $B_n$ has faithful representations by automorphisms of the free group $F_n$ of rank $n$. Let us denote the image of such a representation by $C_n$.
\medskip

$\bullet$ Every group $C_n$ has a homomorphism to the group of $n \times n$ matrices over one-variable Laurent polynomials. This homomorphism is not necessarily injective. 
\medskip

Matrices mentioned in the last bullet point are obtained as follows. Let $\{x_1, \ldots, x_n\}$ be generators of the ambient free group $F_n$. First one computes the $n \times n$ Jacobian matrix $J_\varphi$ of a given automorphism $\varphi \in Aut(F_n)$. This is a matrix of partial Fox derivatives $d_i(y_j)$, where $y_j=\varphi(x_j)$. Fox derivatives are elements of the group ring $\Z F_n$, see Section \ref{Fox} for more details.

The matrix $J_\varphi$ has some interesting properties similar to the ``usual" Jacobian matrix of a multivariate function; for example, given a homomorphism  $\varphi: F_n \to F_n$, $J_\varphi$ is invertible if and only if $\varphi$ is invertible, i.e., is an automorphism of $F_n$ \cite{BJ}. Also, the rows of $J_\varphi$ are linearly independent over the group ring $\Z F_n$ if and only if $\varphi$ is injective \cite{Sh}.

However, the map $\varphi \to J_\varphi$ is not a homomorphism since by the chain rule for Fox derivatives, one has $J_{\varphi \psi} = J_\varphi^\psi J_\psi$, where $J_\varphi^\psi$ denotes the result of applying the homomorphism $\psi$ to all entries of $J_\varphi$. (Any homomorphism of $F_n$ can be extended to the group ring $\Z F_n$ by linearity.)

If we want to get a representation of automorphisms from $C_n$ using Jacobian matrices, we need to apply a homomorphism, call it $\alpha$, to the product $J_\varphi^\psi J_\psi$ to get
$(J_\varphi^\psi J_\psi)^\alpha = J_\varphi^{\psi \alpha} J_\psi^\alpha$. Then, if $J_\varphi^{\psi \alpha} = J_\varphi^\alpha$ for all $\varphi, \psi \in C_n$, we get a
representation of automorphisms from $C_n$ by taking $\varphi \in C_n$ to $J_\varphi^{\alpha}$.

For some $C_n$, this $\alpha$ can be the homomorphism from the group $F_n$ to the infinite cyclic group $<t>$ obtained by taking every $x_i$ to $t$. This homomorphism can be naturally extended to the homomorphism from the group ring $\Z F_n$ to the group ring $\Z [t]$; the latter is the ring of Laurent polynomials over $\Z$.

For a particular representation of $B_n$ by a subgroup of $Aut(F_n)$, known as the Artin representation, the representation by matrices over Laurent polynomials described above is known as the {\it Burau representation}, see \cite{Burau} or \cite{Birman}. It happens so that the g.c.d. of all minors of the same {\it corank} of the [Burau matrix of a braid minus the identity matrix] are invariant under Markov moves and therefore are isotopic invariants of the corresponding links. For (nonzero) minors of the maximum rank, these invariants are known as Alexander polynomials, although the original way of defining Alexander polynomials of a knot was different; it was based on the Wirtinger presentation of the fundamental group $G$ of a knot, see e.g. \cite{Crowell}. Thus, Alexander polynomials are actually invariants of (the isomorphism class of) the group $G$ and therefore cannot possibly distinguish two knots with isomorphic fundamental groups. However, we argue in Section \ref{how} that since, informally speaking, Markov moves form a relatively small subset of the set of all Tietze transformations, our approach in Section \ref{Invariance} allows for a more delicate analysis of how a Burau matrix is affected by Markov moves.

More recently, Wada \cite{W} discovered several other representations of the braid group $B_n$ by automorphisms of $F_n$. These representations were later shown to be faithful  \cite{Shpil} (although this does not play a role in the present paper). Based on Wada's representations, one can obtain other representations of braids by $n \times n$ matrices over Laurent polynomials and produce the corresponding isotopic invariants of knots and links. However, as we conjecture in Section \ref{examples}, we do not get brand new invariants that way; what we get is most likely a specialization of Alexander polynomials.

\section{Preliminaries}\label{Preliminaries}

All facts in this section are well known and can be found, for example, in \cite{Birman}, but we give a concise exposition here for the reader's convenience.

\subsection{Fox derivatives} \label{Fox}

Let $F_n$ be the free group of rank $n$  and $\{x_1, \ldots, x_n\}$ a fixed set of generators. Let $\Z F_n$ be the integral group ring of the group $F_n$.

Partial Fox derivatives $\partial_i$ can be defined for elements of the group $F_n$ using the following  rules, and then extended by linearity to the whole group ring $\Z F_n$: (1) $\partial_i(x_j)=\delta_{ij}$ (Kronecker's delta); ~(2) if $u=vx_i \in F_n$, then $\partial_i(u) = v+ \partial_i(v)$; ~(3) if $u=vx_i^{-1} \in F_n$, then $\partial_i(u) = - vx_i^{-1} + \partial_i(v)$; ~(4) if $u=vx_j, ~j \ne i$, then $\partial_i(u) = \partial_i(v)$.
For example, if $u=x_1 x_2 x_1^{-1} x_2$, then $\partial_1(u) = -x_1 x_2 x_1^{-1} + 1$.

\subsection{Markov moves} \label{Markov}

We will denote braids and the corresponding braid words (i.e., words in the standard  generators $\sigma_i$ of a braid group $B_n$) by the same letters when there is no confusion.

Markov's theorem (see e.g. \cite{Birman}) is:
two braids $\beta_1 \in B_n$ and $\beta_2 \in B_m$ produce isotopically equivalent links when closed if and only if the corresponding braid word $\beta_1$ can be taken to the braid word $\beta_2$ by a sequence of Markov moves, and the latter are:
\medskip

\noindent {\bf 1.} Conjugation in a braid group. That is, if $\beta \in B_n$, one can replace $\beta$ by $\gamma^{-1} \beta \gamma$ for some $\gamma \in B_n$.
\medskip

\noindent {\bf 2.} ``Stabilization". That is, if $\beta \in B_n$, one can multiply $\beta$ by $\sigma_n$ or by $\sigma_n^{-1}$ on the right. Note that the group $B_{n+1}$ is generated by $\sigma_1, \ldots, \sigma_n$, whereas the group $B_{n}$ is generated by $\sigma_1, \ldots, \sigma_{n-1}$.
\medskip

\noindent {\bf 3.} Converse of (2). That is, if it happens so that $\beta = \gamma \sigma_n^{\pm 1}$, where $\gamma \in B_n$, then one can replace $\beta$ by $\gamma$.

\subsection{Matrices over $\Z[t^{\pm 1}]$ and their elementary ideals}\label{minors}

In knot theory, there is a well-known construction of the {\it Alexander matrix} from the Wirtinger presentation of the fundamental group of a knot by generators and defining relations, see e.g. \cite[Chapter 7]{Crowell}. It is similar to our construction of the Burau representation described in the Introduction, in the sense that the Alexander matrix, too, is a matrix of abelianized partial Fox derivatives, but these derivatives are of the defining relations in the Wirtinger presentation.

Since any two presentations (by generators and defining relations) of the same group are equivalent under Tietze transformations, one can study the effect of Tietze transformations on properties of the Alexander matrix and derive knot invariants that way. 

Our approach here is similar, except that we deal here with different matrices over $\Z[t^{\pm 1}]$ (in particular, our matrices are always square), and invariance that we want is under Markov moves, not under Tietze transformations. However, the linear algebra part of our approach is very similar to \cite[Chapter 7.4]{Crowell}.

Specifically, let $M$ be an $n\times n$ matrix over $\Z[t^{\pm 1}]$ and for $0 < n-k <n$, let $E_k=E_k(M)$ be the ideal of the ring $\Z[t^{\pm 1}]$ generated by all minors of size $(n-k)$ of the matrix $M$. Additionally, let $E_k(M) = \{0\}$ if $k<0$ and $E_k(M) = \Z[t^{\pm 1}]$ if $k \ge n$. Then the chain of the ideals $E_0 \subseteq E_1 \subseteq \ldots \subseteq E_n =  E_{n+1} = \Z[t^{\pm 1}]$ is invariant under the usual elementary operations on the rows and/or columns of $M$, as well as under augmenting $M$ by simultaneously adding an extra row $(0 \ldots 0 1)$ (at the bottom) and an extra column $(0 \ldots 0 1)$ (on the right), thus increasing the size of $M$ by 1.

\section{Artin's representation of braid groups and the corresponding polynomial invariants}\label{Artin}

There is a well-known representation, due to Artin, of the braid group $B_n$ in the
group  $Aut(F_n)$ of automorphisms of the free group $F_n$ (see e.g.
\cite[p.25]{Birman}).

Let $F_n$ be generated by $x_1, \ldots, x_n$. Then the automorphism $\hat \sigma_i$
corresponding to the braid generator $\sigma_i$ takes  $x_i$ to
$x_i x_{i+1} x_i^{-1}$,  ~$x_{i+1}$ to $x_i$, and fixes all other generators $x_k$. Denote this representation by $\varphi$.

The corresponding Jacobian matrix $J_\varphi$ is ``mostly" the  identity matrix, with the exception of a $2 \times 2$ cell whose main diagonal is part of the main diagonal of $J_\varphi$, and this cell looks like this:

$$\left(
 \begin{array}{cc} 1-x_i x_{i+1}x_i^{-1} & x_i \\ 1 & 0 \end{array} \right).$$

\noindent If we now define the homomorphism $\alpha$ by taking each $x_i$ to $t$, we will have the condition $J_\varphi^{\psi \alpha} = J_\varphi^\alpha$ satisfied for all $\varphi, \psi$ in the image of Artin's representation, and therefore $\varphi \to J_\varphi^{\alpha}$ will be a homomorphism. The corresponding $2 \times 2$ cell in $J_\varphi^{\alpha}$ will look like this: $\left(
 \begin{array}{cc} 1-t & t \\ 1 & 0 \end{array} \right).$ We note that the inverse of this cell is $\left( \begin{array}{cc} 0 & 1 \\ t^{-1} & 1-t^{-1} \end{array} \right).$

This particular representation of braids by matrices $J_\varphi^{\alpha}$ is known as the {\it Burau representation}, see \cite{Burau} or \cite{Birman}.

\subsection{Invariance under Markov moves}\label{Invariance}

Now we are going to show that the ideals of the ring $\Z[t^{\pm 1}]$ generated by all $k \times k$ minors of the matrix $J_\varphi^\alpha - I$ are invariant under Markov moves applied to the braid that corresponds to the automorphism $\varphi$. Here $I$ denotes the identity matrix of the right size.

Invariance under conjugation is well known, so we are going to study the effect (on the matrix $J_\varphi^\alpha - I$) of multiplying a braid $\beta$ from $B_{n}$ by $\sigma_n$ or by $\sigma_n^{-1}$ on the right. Denote $J_\varphi^\alpha$ by $J_\beta$ to simplify the notation.


	\begin{align*}
		J_{\beta \sigma_n}- I&=\begin{pmatrix}
			a_{11}	&\dots	&	&\vdots\\
					&\ddots	&y	&0\\
					&\dots	&x	&0\\
					&\dots	&0	&1
		\end{pmatrix}		\begin{pmatrix}
			1	&\dots	&	&\vdots\\
					&\ddots	&0	&0\\
					&\dots	&1-t	&t\\
					&\dots	&1 & 0\\
		\end{pmatrix} - I\\
		&=\begin{pmatrix}
			a_{11}-1	&\dots	&	&\vdots\\
							&\ddots	&y-yt	&yt\\
							&\dots	&x-xt-1	&xt\\
				0	&\dots	&1	&-1\\
		\end{pmatrix}.\\
	\end{align*}

After adding the last column to the second column from the right, we get the matrix

$$\begin{pmatrix}
			a_{11}-1	&\dots	&	&\vdots\\
							&\ddots	&y	&yt\\
							&\dots	&x-1	&xt\\
				0	&\dots	&0	&-1\\
		\end{pmatrix}.$$
\medskip

Thus, $\det(J_{\beta \sigma_n} - I) = -det(J_\beta - I)$, so $\det(J_\beta - I)$ and $\det(J_{\beta \sigma_n} - I)$ generate the same ideal of $\Z[t^{\pm 1}]$. 
Also,  chains of the elementary ideals of the matrices $(J_{\beta}-I)$ and $(J_{\beta \sigma_n}-I)$ are the same since the latter matrix can be obtained from the former by a sequence of elementary operations, see Section \ref{minors}.

Now we compare the matrices $J_\beta - I$ and $J_{\beta \sigma_n^{-1}} - I$.

	\begin{align*}
		 J_{\beta \sigma_n^{-1}} - I = \begin{pmatrix}
			a_{11}	&\dots	&	&\vdots\\
					&\ddots	&y	&0\\
					&\dots	&x	&0\\
					&\dots	&0	&1
		\end{pmatrix}
		\begin{pmatrix}
			1	&\dots	&	&\vdots\\
					&\ddots	&0	&0\\
					&\dots	&0	&1\\
					&\dots	&t^{-1} & 1-t^{-1}\\
		\end{pmatrix} - I
		&=\begin{pmatrix}
			a_{11}	&\dots	&	&\vdots\\
					&\ddots	&0	&y\\
					&\dots	&0	&x\\
			0		&\dots	&t^{-1} & 1-t^{-1}
		\end{pmatrix} - I\\ =
\begin{pmatrix}
			a_{11}-1	&\dots	&	&\vdots\\
					&\ddots	&0	&y\\
					&\dots	&-1	&x\\
			0		&\dots	&t^{-1} & -t^{-1}
\end{pmatrix}.
	\end{align*}

After adding the second column from the right to the rightmost column and then switching the last two columns, we get the matrix

         $\begin{pmatrix}
			a_{11}-1	&\dots	&	&\vdots\\
					&\ddots	&y	&0\\
					&\dots	&x-1	&-1\\
			0		&\dots	&0 & -t^{-1}
        \end{pmatrix}$.  We see that $\det(J_{\beta \sigma_n^{-1}}- I) = t^{-1}\det(J_{\beta}-I)$. Also,  chains of the elementary ideals of the matrices $(J_{\beta}-I)$ and $(J_{\beta \sigma_n^{-1}}-I)$ are the same since the latter matrix can be obtained from the former by a sequence of elementary operations, see Section \ref{minors}.

\section{Can Burau matrices be used to show that
the right and left trefoil knots are not isotopic?}\label{how}

It is a common belief that Alexander matrices, being obtained from presentations of the fundamental group of a knot, cannot be used to distinguish two knots with isomorphic
fundamental groups. However, ``invariance under isomorphisms" is typically established (see e.g. \cite[Chapter 7.4]{Crowell}) as invariance under Tietze transformations.

On the other hand, Markov moves, informally speaking, form a relatively small subset of the set of all Tietze transformations, and this is why our approach in Section \ref{Invariance} allows for a more delicate analysis.

Specifically, let us illustrate our point using the example of the right and left trefoil knots.

The braid corresponding to the right trefoil knot is $\beta=\sigma_1^3$,  and the corresponding matrix $J_{\beta} - I$ is

$$\left(
 \begin{array}{cc} -t^3+t^2-t & t^3-t^2+t \\ t^2-t+1 & -t^2+t-1 \end{array} \right).$$

On the other hand, the braid corresponding to the left trefoil knot is $\beta'=\sigma_1^{-3}$,  and the corresponding matrix $J_{\beta'} - I$ is

$$\left(
 \begin{array}{cc} -t^{-2}+t^{-1}-1 & t^{-2}-t^{-1}+1 \\ t^{-3}-t^{-2}+t^{-1}  & -t^{-3}+t^{-2}-t^{-1} \end{array} \right).$$

It seems plausible that no sequence of Markov moves applied to the braid $\beta'$ can entirely eliminate monomials $t^k$ with $k<0$ from the matrix $J_{\beta'} - I$. It is clear from our Section \ref{Invariance} that this is true for Markov moves of type 2 and 3 (``stabilization" and its converse), but the effect of the conjugation is more elusive. Note that conjugation can be done not by just any matrix over $\Z[t^{\pm 1}]$, but only by Burau matrices of braids.
%
%

That said, we know that the braids $\sigma_1$ and  $\sigma_1^{-1}$ produce the same knot (the unknot), and the braids $\sigma_1^{2}$ and  $\sigma_1^{-2}$ produce the same link (the Hopf link). All the entries of the Burau matrix of $\sigma_1^{-2}$ have monomials $t^k$ with $k<0$, and yet there is a sequence of Markov moves that takes $\sigma_1^{-2}$ to $\sigma_1^{2}$, even though the Burau matrix of $\sigma_1^{2}$ does not have any monomials $t^k$ with $k<0$. It may therefore be useful to find an explicit sequence of Markov moves that takes $\sigma_1^{-2}$ to $\sigma_1^{2}$ to try to find an algebraic reason why there is no such sequence taking $\sigma_1^{-3}$ to $\sigma_1^{3}$.

\subsection{A group associated to an endomorphism of $F_n$}\label{torus}

We mentioned above that ``Markov moves form a subset of the set of all Tietze transformations", but Tietze transformations are applied to a group presentation by generators and defining relations, whereas Markov moves are applied to elements of a braid group $B_n$. Thus, an explanation is in order.

Let $F_n$ be the free group of rank $n$ with a set $\{x_1, \ldots, x_n\}$ of free generators.
Let $\varphi: F_n \to F_n$ be an endomorphism of $F_n$. Let $\varphi(x_i)=y_i$. Then we can associate the following group, given by generators and defining relations, to the endomorphism $\varphi$:

$$G_\varphi = \langle x_1, \ldots, x_n, y_1, \ldots, y_n ~|~ y_1=x_1, \ldots, y_n=x_n \rangle.$$

\noindent Equivalently, $G_\varphi = \langle x_1, \ldots, x_n, y_1, \ldots, y_n ~|~ x_1^{-1} y_1, \ldots, x_n^{-1} y_n \rangle.$ This group does not have a special name, to the best of our knowledge. It should not be confused with the {\it mapping torus} of $\varphi$, see e.g. \cite{Mutanguha}.

The Alexander matrix $A$ (see e.g. \cite[Chapter 7.3]{Crowell}) of the latter presentation of $G_\varphi$ can be obtained from the matrix $J_\varphi$ by multiplying each entry in row $i$ by $x_i^{-1}$ (on the left) and then subtracting the diagonal matrix with $x_i^{-1}$ in the $(i, i)$th position. Thus, after applying the abelianization map $\alpha: x_i \to t$, we get  $A^\alpha = t^{-1} (J^\alpha_\varphi - I)$. Therefore, the chain of elementary ideals of the matrix $A^\alpha$ is the same as that of the matrix $(J^\alpha_\varphi - I)$.

Tietze transformations applied to the presentation of the group $G_\varphi$ yield elementary operations on rows and columns of the matrix $A^\alpha$, just as Markov moves yield elementary operations on rows and columns of the matrix $(J^\alpha_\varphi - I)$. However, since Markov moves of type (1) are just conjugations, not arbitrary isomorphisms, their effect on the matrix $(J^\alpha_\varphi - I)$ can be controlled better (at least in theory).

\section{Wada's representation of braid groups and the corresponding knot invariants}
\label{Wada}

Wada \cite{W} found several other representations of the braid group $B_n$ in the
group  $Aut(F_n)$; these were later shown to be faithful \cite{Shpil}.

Of interest to us in this paper is the following Wada's representation. In this  representation, the automorphism $\hat \sigma_i$
corresponding to the braid generator $\sigma_i$, takes  $x_i$ to
$x_i^2 x_{i+1}$,  ~$x_{i+1}$ to $x_{i+1}^{-1} x_i^{-1} x_{i+1}$, and fixes all other generators $x_k$.

Again, the corresponding Jacobian matrix $J_{\hat \sigma_i}$ is ``mostly" the $n \times n$ identity matrix, with the exception of a $2 \times 2$ cell whose main diagonal is part of the main diagonal of $J_\varphi$, and this cell looks like this:

$$\left(
 \begin{array}{cc} 1+x_i & x_i^2 \\ -x_{i+1}^{-1} x_i^{-1} & -x_{i+1}^{-1} + x_{i+1}^{-1} x_i^{-1} \end{array} \right).$$

\noindent Here defining the abelianization homomorphism $\alpha$ to satisfy the condition $J_\varphi^{\psi \alpha} = J_\varphi^\alpha$ for all $\varphi, \psi$ in the image of Wada's representation is a little more tricky. It can be defined as follows. If $i \ge 1$ is odd, then $\alpha(x_i)=t$; otherwise, $\alpha(x_i)=t^{-1}$. Then the corresponding $2 \times 2$ cell in $J_\varphi^{\alpha}$ will look like this if $i$ is odd:
\medskip

$\left(
 \begin{array}{cc} 1+t & t^2 \\ -1 & 1-t \end{array} \right).$ The inverse of this cell is
 $\left(
 \begin{array}{cc} 1-t & -t^2 \\ 1 & 1+t \end{array} \right).$
\medskip

\noindent If $i$ is even, the cell will look like this:
\medskip

$\left(
 \begin{array}{cc} 1+t^{-1} & t^{-2} \\ -1 & 1-t^{-1} \end{array} \right).$ The inverse of this cell is $\left(
 \begin{array}{cc} 1-t^{-1} & -t^{-2} \\ 1 & 1+t^{-1}  \end{array} \right).$

\medskip

Denote the matrix $J_\varphi^{\alpha}$ corresponding to a braid $\beta$ by just $J_\beta$.
We are now going to see the effect of Markov moves applied to the braid $\beta$ on the matrix $(J_\beta - I)$. First, let us assume that $\beta$ is multiplied by $\sigma_n$ with an odd $n$.

	\begin{align*}
		 J_{\beta \sigma_n} - I = \begin{pmatrix}
			a_{11}	&\dots	&	&\vdots\\
					&\ddots	&y	&0\\
					&\dots	&x	&0\\
					&\dots	&0	&1
		\end{pmatrix}
		\begin{pmatrix}
			1	&\dots	&	&\vdots\\
					&\ddots	&0	&0\\
					&\dots	&1+t	&t^2\\
					&\dots	&-1 & 1-t\\
		\end{pmatrix} - I
		&=\begin{pmatrix}
			a_{11}	&\dots	&	&\vdots\\
					&\ddots	&y(1+t)	&yt^2\\
					&\dots	&x(1+t)	&xt^2\\
			0		&\dots	&-1 & 1-t
		\end{pmatrix} - I\\ =
\begin{pmatrix}
			a_{11}-1	&\dots	&	&\vdots\\
					&\ddots	&y(1+t)	&yt^2\\
					&\dots	&x(1+t)-1	&xt^2\\
			0		&\dots	&-1 & -t
\end{pmatrix}.
	\end{align*}

Now we add the last row multiplied by $xt$ to the second row from the bottom, and then add the last row multiplied by $yt$ to the third row from the bottom. This gives $\begin{pmatrix}
			a_{11}-1	&\dots	&	&\vdots\\
					&\ddots	&y	&0\\
					&\dots	&x-1	&0\\
			0		&\dots	&-1 & -t
\end{pmatrix}$. Finally, we multiply the rightmost column by $-t^{-1}$ and add it to the second
column from the right to get $\begin{pmatrix}
			a_{11}-1	&\dots	&	&\vdots\\
					&\ddots	&y	&0\\
					&\dots	&x-1	&0\\
			0		&\dots	&0 & -t
\end{pmatrix}$.
\medskip

We see that $\det(J_{\beta \sigma_n}- I) = -t\det(J_{\beta}-I)$. Also,  chains of the elementary ideals of the matrices $(J_{\beta}-I)$ and $(J_{\beta \sigma_n}-I)$ are the same since the latter matrix can be obtained from the former by a sequence of elementary operations, see Section \ref{minors}.
\medskip

Now we compare the matrices $J_\beta - I$ and $J_{\beta \sigma_n^{-1}} - I$.

	\begin{align*}
		 J_{\beta \sigma_n^{-1}} - I = \begin{pmatrix}
			a_{11}	&\dots	&	&\vdots\\
					&\ddots	&y	&0\\
					&\dots	&x	&0\\
					&\dots	&0	&1
		\end{pmatrix}
		\begin{pmatrix}
			1	&\dots	&	&\vdots\\
					&\ddots	&0	&0\\
					&\dots	&1-t	&-t^2\\
					&\dots	&1 & 1+t\\
		\end{pmatrix} - I
		&=\begin{pmatrix}
			a_{11}	&\dots	&	&\vdots\\
					&\ddots	&y(1-t)	&-yt^2\\
					&\dots	&x(1-t)	&-xt^2\\
			0		&\dots	&1 & 1+t
		\end{pmatrix} - I\\ =
\begin{pmatrix}
			a_{11}-1	&\dots	&	&\vdots\\
					&\ddots	&y(1-t)	&-yt^2\\
					&\dots	&x(1-t)-1	&-xt^2\\
			0		&\dots	&1 & t
\end{pmatrix}.
	\end{align*}
\medskip

Now we add the last row multiplied by $xt$ to the second row from the bottom, and then add the last row multiplied by $yt$ to the third row from the bottom. This gives $\begin{pmatrix}
			a_{11}-1	&\dots	&	&\vdots\\
					&\ddots	&y	&0\\
					&\dots	&x-1	&0\\
			0		&\dots	&1 & t
\end{pmatrix}$.  Finally, we multiply the rightmost column by $-t^{-1}$ and add it to the second
column from the right to get $\begin{pmatrix}
			a_{11}-1	&\dots	&	&\vdots\\
					&\ddots	&y	&0\\
					&\dots	&x-1	&0\\
			0		&\dots	&0 & t
\end{pmatrix}$.
\medskip

We see that $\det(J_{\beta \sigma_n^{-1}}- I) = t\det(J_{\beta}-I)$. Also,  chains of the elementary ideals of the matrices $(J_{\beta}-I)$ and $(J_{\beta \sigma_n^{-1}}-I)$ are the same since the latter matrix can be obtained from the former by a sequence of elementary operations, see Section \ref{minors}.

This completes the case where $n$ is odd. The case where $n$ is even is treated similarly since it amounts to just replacing $t$ by $t^{-1}$ in the above.

\subsection{Wada's polynomials}\label{polynomials}
Invariants coming from Wada's representation are ideals $E_k$ of the Laurent polynomial ring $\Z[t^{\pm 1}]$ generated by all  $(n-k) \times (n-k)$ minors of the relevant $n \times n$ matrix $J_\varphi^{\alpha} - I$.  If the $n \times n$ minor (i.e., the determinant) of the matrix $J_\varphi^{\alpha} - I$ is 0, then we consider $(n-1) \times (n-1)$  minors, etc.

The polynomial that generates nonzero $E_k$ with the smallest $k \ge 0$ is what we call the Wada polynomial of the relevant knot or link.
To avoid ambiguity, we normalize Wada polynomials as follows: (1) the polynomial should not contain any negative exponents on $t$; (2) the constant term should not be 0; (3) the coefficient at the highest degree monomial should be positive. For example, when we normalize the polynomial $t^{-2}-1-t$, we multiply it by $t^{2}$, change the sign and get $t^3+t^2-1$.

\section{Examples}\label{examples}

Here we give examples of Wada polynomials of some simple knots and links. We also compare them to Alexander polynomials.

\subsection{The unknot} The braid that corresponds to the unknot is $\beta=\sigma_1$,  and the corresponding matrix $J_{\beta} - I$ is

$$\left(
 \begin{array}{cc} t & t^2 \\ -1 & -t \end{array} \right).$$

The determinant of this matrix is 0, so we look at the g.c.d. of $1 \times 1$  minors, and this is equal to 1, which is the same as the Alexander polynomial of the unknot.

\subsection{The Hopf link} The braid that corresponds to the Hopf link is $\beta=\sigma_1^2$,  and the corresponding matrix $J_{\beta} - I$ is

$$\left(
 \begin{array}{cc} 2t & 2t^2 \\ -2 & -2t \end{array} \right).$$

The determinant of this matrix is 0, so we look at the g.c.d. of $1 \times 1$ minors, and this is equal to 2. This is the Wada polynomial of the Hopf link. Note that the Alexander polynomial of the Hopf link is $1-t$.

\subsection{The trefoil knot} The braid that corresponds to the (right) trefoil knot is $\beta=\sigma_1^3$,  and the corresponding matrix $J_{\beta} - I$ is

$$\left(
 \begin{array}{cc} 3t & 3t^2 \\ -3 & -3t \end{array} \right).$$

Again, we compute the g.c.d. of $1 \times 1$  minors, and this is equal to 3. This is the Wada polynomial of the trefoil knot. Note that the Alexander polynomial of the trefoil knot is $1-t+t^2$.

\subsection{$(2, k)$ torus knots} The computation of the matrix $J_{\beta}$ easily generalizes by induction to knots (or links) corresponding to the braids $\sigma_1^k$. If $k \ge 3$ is odd, this gives $(2, k)$ torus knots. If $k \ge 2$ is even, we have $(2, k)$ torus links. In either case,
the Wada polynomial is equal to $k$. The Alexander polynomial is $\frac{1+t^k}{1+t} = 1-t+t^2- \ldots + (-1)^{k-1}t^{k-1}$.

\subsection{The figure eight knot}  The braid that corresponds to the figure eight knot is $\beta=(\sigma_1 \sigma_2^{-1})^2$,  and the corresponding matrix $J_{\beta} - I$ is

$$\left(
 \begin{array}{ccc} 3t & 3t^2-4t & -4 \\ -3+t^{-1} & 5-3t-3t^{-1} &  -3t^{-2}+ 4t^{-1} \\
 -1 & 3-t & 3t^{-1} \end{array} \right).$$

The determinant of this matrix is 0, and the g.c.d. of all $2 \times 2$ minors is 5, so this is the Wada polynomial of the figure eight knot.  The Alexander polynomial of the figure eight knot is $t^{2}-3t+1$.

\subsection{The square knot} The square knot is a composition of two copies of the right trefoil knot. The corresponding braid is $\beta=\sigma_1^3 \sigma_2^3$, and the corresponding matrix $J_{\beta} - I$ is

$$\left(
 \begin{array}{ccc} 3t & 9t+ 3t^2 & 9 \\ -3 & 3t^{-1}-3t-9 & 3t^{-2}- 9t^{-1} \\
 0 & -3 & -3t^{-1} \end{array} \right).$$

The determinant of this matrix is 0, and the g.c.d. of all $2 \times 2$ minors is 9, so this is the Wada polynomial of the square knot. The Alexander polynomial of the square knot is
$(1-t+t^2)^2$.

\subsection{Granny's knot} Granny's knot is a composition of the right and left trefoil knots. The corresponding braid is $\beta=\sigma_1^3 \sigma_2^{-3}$, and the corresponding matrix $J_{\beta} - I$ is

$$\left(
 \begin{array}{ccc} 3t & -9t+ 3t^2 & -9 \\ -3 & -3t^{-1}-3t+9 & -3t^{-2}+ 9t^{-1} \\
 0 & 3 & 3t^{-1} \end{array} \right).$$

The determinant of this matrix is 0, and the g.c.d. of all $2 \times 2$ minors is 9, so this is the Wada polynomial of granny's knot, the same as that of the square knot. The Alexander polynomial of granny's knot is $(1-t+t^2)^2$, the same as that of the square knot.

\medskip

The above examples make it appear likely that the following conjecture holds:\footnote{Tetsuya Ito has pointed out to us that a proof of this conjecture can be recovered from results of \cite{W} and \cite{Sakuma}.}

\begin{conjecture}
The Wada polynomial is the specialization of the Alexander polynomial at $t=-1$.
\end{conjecture}

\section{Supplement: a representation by matrices over $\Z[s^{\pm 1}, t^{\pm 1}]$}\label{Supplement}

In this section, we point out a representation of braid groups $B_n$ by $n \times n$ matrices that is not in line with the main theme of the present paper (since it does not come from representing braid groups by automorphisms). Like other representations by $n \times n$  matrices considered in this paper, it is ``local" in the sense that each braid generator $\sigma_i$ is represented by an  $n \times n$ matrix that differs from the identity matrix just by a $2 \times 2$ cell in the right place along the main diagonal.

We note, in passing, that although braid groups are known to be linear \cite{Bigelow01}, \cite{Krammer}, it is still unknown whether or not they have a faithful representation by matrices over $\Q$. If they do, this would imply, in particular, that the word problem in the group $B_n$ is solvable in quasilinear time, see \cite{Olsh}.

The representation we are talking about appears in \cite{Pourkia} (also see references therein). In this case, the $2 \times 2$ cell mentioned above is:


$$\left(
 \begin{array}{cc} 1-st & t \\ s & 0 \end{array} \right).$$

\noindent The Burau representation is the specialization of this representation at $s=1$.

Denote by $M_\beta$ the matrix (over $\Z[s^{\pm 1}, t^{\pm 1}]$) that corresponds to the braid $\beta$ under this representation.

It can be shown the same way this was done in Section \ref{Invariance} that the chain of the elementary ideals (of the ring $\Z[s^{\pm 1}, t^{\pm 1}]$) of the matrix $(M_\beta - I)$ is invariant under Markov moves, and therefore one can get the corresponding knot/link invariants this way. The ring $\Z[s^{\pm 1}, t^{\pm 1}]$ is not a principal ideal domain, so there are ideals that are not generated by a single polynomial. However, for ideals that correspond to braid words under the above representation this seems to be the case. Moreover, if one denotes the Alexander polynomial that corresponds to a braid $\beta$ by $AL_\beta(t)$, then the ``leading" invariant corresponding to our representation of the same braid $\beta$ seems to be $AL_\beta(st)$. For example, the ``leading" invariant of the Hopf link will be the ideal of $\Z[s^{\pm 1}, t^{\pm 1}]$ generated by the polynomial $1-st$, whereas the ``leading" invariant of the trefoil knot (both left and right) will be the ideal generated by the polynomial $1-st+s^2t^2$.


Therefore, the above representation does not seem to yield new invariants of knots/links. However, it may be of interest for a different reason. The Burau representation of the braid group $B_n$ is known not to be faithful if $n \ge 5$ \cite{Bigelow99}, \cite{Long} and faithful if $n=3$ \cite{MagnusPeluso}. The kernel of the representation in this section cannot be larger than that of the Burau representation since the Burau representation is a specialization of it. Whether or not it is strictly smaller (or even trivial) is an interesting question.



\medskip

\noindent  {\bf Acknowledgement.} 
I am grateful to Lorenzo Traldi for helpful comments.


\baselineskip 11 pt


\begin{thebibliography}{ABC}


\bibitem{Bigelow99}
S. Bigelow, {\it The Burau representation is not faithful for $n=5$}, Geom. Topol. {\bf 3} (1999), 397--404 (electronic).

\bibitem{Bigelow01}
S. Bigelow, {\it Braid groups are linear}, J. Amer. Math. Soc. {\bf 14} (2001),
471--486.

\bibitem{BJ}
J.~S. Birman, {\it  An~inverse function theorem for free groups}, \emph{Proc.\ Amer.\
  Math.\ Soc.{}}, \textbf{41} (1973), 634--638.

\bibitem{Birman}
J. S. Birman, {\sl Braids, links and  mapping  class  groups},
Ann.  Math. Studies {\bf 82}, Princeton Univ. Press, 1974.

\bibitem{Burau}
W. Burau, {\it \"Uber Zopfgruppen und gleichsinnig verdrillte Verkettungen}, Abh. Math. Sem. Univ. Hamburg {\bf 11} (1936), 179--186.

\bibitem{Crowell}
R. H. Crowell, R. H. Fox, {\sl Introduction to knot theory}, Grad. Texts in Math., No. {\bf 57}.  Springer-Verlag, New York-Heidelberg, 1977.

\bibitem{Fox}
R.~H. Fox, {\it Free differential calculus. {I.} {D}erivation in the~free group ring}, Ann.\ Math.{}~(2), \textbf{57} (1953), 547--560.

\bibitem{Krammer}
D. Krammer, {\it Braid groups are linear}, Ann. of Math. (2) {\bf 155} (2002),  131--156.

\bibitem{Long}
D. D. Long, M. Paton, {\it The Burau representation is not faithful for $n \ge 6$},  Topology {\bf 32} (1993), 439--447.

\bibitem{MagnusPeluso}
W. Magnus, A. Peluso, {\it On a theorem of V. I. Arnold},
Comm. Pure Appl. Math. {\bf 22} (1969), 683--692.

\bibitem{Mutanguha}
J. P. Mutanguha, {\it Irreducibility of a free group endomorphism is a mapping torus invariant}, Comment. Math. Helv. {\bf 96} (2021), 47--63.

\bibitem{Olsh}
A. Olshanskii, V. Shpilrain, {\it Linear average-case complexity of algorithmic problems in groups,} J. Algebra {\bf 668} (2025), 390--419.

\bibitem{Pourkia}
A. Pourkia, {\it A basic $n$-dimensional representation of Artin braid group $B_n$, and a general Burau representation}, Mathematics and Statistics {\bf 10} (2022), 145--152.

\bibitem{Sakuma}
M. Sakuma, {\it A note on Wada's group invariants of links},  Proc. Japan Acad. Ser. A Math. Sci. {\bf 67} (1991), 176--177.

\bibitem{Sh}
V. Shpilrain, {\it On monomorphisms of free groups}, Arch. Math. {\bf 64} (1995), 465--470.

\bibitem{Shpil}
V. Shpilrain, {\it Representing braids
by automorphisms}, Internat. J. Algebra and Comput. {\bf  11} (2001),
773--778.


\bibitem{W}
M. Wada,  {\it Group invariants of links},  Topology {\bf 31} (1992),  399--406.



\end{thebibliography}
\end{document}